\documentclass[a4paper,11pt]{article}
\usepackage[super,square]{natbib}
\usepackage{anysize}
\marginsize{2cm}{2cm}{1cm}{2cm}

\usepackage{fancyhdr}

\usepackage{soul}
\usepackage{graphicx}
\usepackage{subcaption}
\usepackage[dvipsnames]{xcolor}
\usepackage{hyperref}
\hypersetup{
    colorlinks=true,
    linkcolor=black,
    citecolor=CadetBlue,
    filecolor=CadetBlue,      
    urlcolor=CadetBlue,
}
\usepackage{lipsum}
\usepackage{multicol}
\usepackage{amsmath, amssymb, amsfonts, physics}

\usepackage{float}
\renewenvironment{abstract}
 {\par\noindent\textbf{\abstractname}\ \ignorespaces \\}
 {\par\noindent\medskip}

\begin{document}
\pagestyle{fancy}
\thispagestyle{empty}
\fancyhead[R]{}
\fancyhead[L]{}
\renewcommand*{\thefootnote}{\fnsymbol{footnote}}
\begin{center}
\Large{\textbf{Stability and Synchronization of Kuramoto Oscillators}}
\vspace{0.4cm}
\normalsize
\\ Abhiram Gorle, Stanford University \\ 
\footnotesize Correspondence at: abhiramg@stanford.edu\\
\vspace{0.1cm}
\medskip
\normalsize
\end{center}
{\color{gray}\hrule}
\vspace{0.4cm}
\begin{abstract}
Imagine a group of oscillators, each endowed with their own rhythm or frequency, be it the ticking of a biological clock, the swing of a pendulum, or the glowing of fireflies. While these individual oscillators may seem independent of one another at first glance, the true magic lies in their ability to influence and synchronize with one another, like a group of fireflies glowing in unison. 

The Kuramoto model was motivated by this phenomenon of collective
synchronization, when a group of a large number of oscillators spontaneously lock to a common frequency, despite vast differences in their individual frequencies (A.T. Winfree 1967, S.H. Strogatz 2000). Inspired by Kuramoto's groundbreaking work in the 1970s, this model captures the essence of how interconnected systems, ranging from biological networks to power grids, can achieve a state of synchronization.

This work aims to study the stability and synchronization of Kuramoto oscillators, starting off with an introduction to Kuramoto Oscillators and it's broader applications. We then at a graph theoretic formulation for the same and establish various criterion for the stability, synchronization of Kuramoto Oscillators. Finally, we broadly analyze and experiment with various physical systems that tend to behave like Kuramoto oscillators followed by further simulations.

\end{abstract}
{\color{gray}\hrule}

\tableofcontents
\section{Introduction}
The Kuramoto model was proposed to study huge populations of coupled limit-cycle oscillators whose natural frequencies are known a priori. 
Tha background work for Kuramoto model is laid in Winfree (1967), where it proposed that \textbf{'each oscillator was coupled to the collective rhythm generated by the whole population'}, analogous to a mean-field approximation in Physics (where we model certain random variables in terms of the mean of their variation). So, Winfree's proposed model for a system of $N$ oscillators was as follows:
\begin{equation*}
\dot{\theta}_i=\omega_i+K\left(\sum_{j=1}^N A_{j i} X\left(\theta_j\right)\right) Z\left(\theta_i\right)
\end{equation*}
where, $K$ is the coupling strength and $A_{j i}$ is a measure of the communication capacity between different channels and $X(\theta_j)$ is the influence of other channels on the $i^{th}$ channel.

But this model was not widely accepted by the scientific community because it lacked certain symmetries like that of translational invariance when both phases are slightly perturbed, the model does not remain invariant. The classical Kuramoto model proposed in 1975 is as follows:

$$
\dot{\theta}_i=\omega_i+\sum_{j=1}^n \Gamma_{i j}\left(\theta_j-\theta_i\right)
$$
Where \(\theta_i\) are the phases and \(\omega_i\) are the limit-cycle frequencies of the oscillators, we use a sine function to couple them. This choice simplifies the model to the following form:
$$
\dot{\theta}_i=\omega_i+\frac{K}{n} \sum_{j=1}^n \sin \left(\theta_j-\theta_i\right)
$$

If the coupling constant \(K\) is sufficiently large, then even an initially incoherent set of oscillators will first gain partial coherence and eventually become fully synchronized. A natural question is: \textbf{How large must the coupling constant be to induce synchronization? And, for any given value of \(K\), how can we measure or quantify the degree of synchronization in the system?}

\subsection{Order parameter}

An important way to measure synchronization in the Kuramoto model is through the \emph{order parameter}. 
Conceptually, this parameter represents the \emph{centroid} of the oscillators 
when each is viewed as a point on the complex unit circle. Mathematically, we write
\[
r\,e^{i\psi} \;=\; \frac{1}{n}\,\sum_{j=1}^{n} e^{\,i\,\theta_j},
\]
where \(\theta_j\) is the phase of the \(j\)-th oscillator, 
and \(r \in [0,1]\) is the \emph{magnitude} of the order parameter. 
A higher value of \(r\) indicates that the oscillators are more closely clustered in phase 
(i.e., more synchronized). Specifically:
\begin{itemize}
    \item When the oscillators are perfectly phase-synchronized, \(r = 1\).
    \item When the oscillators are spread evenly around the unit circle, \(r = 0\).
    \item For intermediate values of \(r \in (0,1)\), the system exhibits partial synchronization, 
    with greater values of \(r\) corresponding to higher levels of cohesiveness.
\end{itemize}

Expressing the original oscillator equations in terms of \(r\) and \(\psi\) can help us analyze 
how the system’s synchronization evolves over time.

Using this formulation, we get:
\begin{equation*}
    \boxed{\dot{\theta}_i=\omega_i+K r \sin \left(\psi-\theta_i\right)}
\end{equation*}

If we assume our critical coupling constant for synchronization to occur to be $K_c$, then the above plot gives us an idea of how order parameter evolves with time.
\begin{figure}
    \centering
    \includegraphics[scale = 1.2]{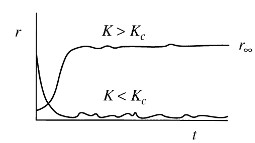}
    \caption{Evolution of $r(t)$}
    \label{fig:orderparam}
\end{figure}
\vspace{8pt}

{\color{gray}\hrule}
\begin{center}
\section{Stability}
\textbf{A primer on Stability of Kuramoto Oscillators}
\end{center}
{\color{gray}\hrule}
\subsection{Mathematical Setup}

Before we analyze stability, let us first set up the problem using a graph-theoretic formulation of the Kuramoto model. 

\vspace{8pt}

\noindent
\textbf{Incidence Matrix and Laplacian:}
Let \(\mathcal{G} = (\mathcal{V},\mathcal{E})\) be a (connected) undirected graph with:
\[
|\mathcal{V}| = N
\quad\text{and}\quad
|\mathcal{E}| = e.
\]
We label each vertex by an index \(i \in \{1,\dots,N\}\) and each edge by an index \(k \in \{1,\dots,e\}\). To work with oriented edges, we fix an arbitrary orientation \(\sigma\), which assigns a direction to each undirected edge.

Consider an oriented graph \(\mathcal{G}^\sigma\) with \(N\) vertices and \(e\) edges. Its incidence matrix \(B\) is the \(N \times e\) matrix such that
\[
B_{ij} \;=\;
\begin{cases}
1 & \text{if edge } j \text{ is incoming to vertex } i,\\
-1 & \text{if edge } j \text{ is outgoing from vertex } i,\\
0 & \text{otherwise}.
\end{cases}
\]
The \emph{Laplacian} of \(\mathcal{G}\), denoted by \(L\), is then defined as
\[
L \;=\; B\,B^T.
\]
Notably, \(L\) is \emph{independent} of the orientation \(\sigma\), and it has several important properties:
\begin{itemize}
    \item  \(L\) is \textit{positive semidefinite} with at least one zero eigenvalue, 
    \item The dimension of its nullspace (i.e., the algebraic multiplicity of eigenvalue 0) is the number of connected components of \(\mathcal{G}\). Since \(\mathcal{G}\) is connected, \(L\) has exactly one zero eigenvalue,
    \item The corresponding eigenvector to the zero eigenvalue is \(\mathbf{1}_N = (1,1,\dots,1)^T\),
    \item The first nonzero eigenvalue \(\lambda_2(L)\) measures the \textit{algebraic connectivity} of \(\mathcal{G}\). 
    \item If a positive weight \(W_k\) is assigned to each edge \(k \in \mathcal{E}\) and constructing \(W = \mathrm{diag}(W_1,\dots,W_e)\), then the \emph{weighted Laplacian} is
\[
L_W(\mathcal{G})
\;=\;
B\,W\,B^T,
\]
which preserves the above properties regarding nullspace, connectivity, etc.
\end{itemize}

\vspace{8pt}

\noindent
\textbf{Dynamics:}
We focus on the following system of differential equations describing the phases \(\theta \in \mathbb{R}^N\):
\begin{equation}
    \dot{\theta} \;=\; \omega \;-\; \frac{K}{N}\,B\,\sin\bigl(B^T \theta\bigr),
\end{equation}
where \(\omega \in \mathbb{R}^N\) is the vector of intrinsic (natural) frequencies, and \(K > 0\) is the coupling strength. 

\vspace{8pt}

\noindent
The order parameter (defined previously) is computed here as follows:
\begin{equation*}
r^2=1-\frac{1}{N}\left[e^{j \theta}\right]^* L\left[e^{j \theta}\right]=1-\frac{1}{N}\left([\cos \theta]^T L[\cos \theta]+[\sin \theta]^T L[\sin \theta]\right)
\end{equation*}
\textbf{Generalized Inverse:}
The (Moore--Penrose) pseudoinverse of \(V^T B\) is denoted by \(\bigl(V^T B\bigr)^{\#}\), and equals
\[
\bigl(V^T B\bigr)^{\#} 
\;=\;
B^T\,V\,\Lambda^{-1},
\]
where \(\Lambda\) is the diagonal matrix of the nonzero eigenvalues of \(L\). Because \(L\) has at least one zero eigenvalue for a connected graph, that eigenvalue is omitted in \(\Lambda\). We also have
\[
\left(\sin \bigl(B^T \theta\bigr)\right)_{R\bigl(B^T\bigr)}
\;=\;
B^T\,V\,\Lambda^{-1}\,V^T\,\frac{N\,\omega}{K}
\;=\;
B^T\,L^{\#}\,\frac{N\,\omega}{K},
\]
where \(L^{\#} = V\,\Lambda^{-1}\,V^T\) is the pseudoinverse of the Laplacian \(L\).
\vspace{8pt}

\noindent
\textbf{The Kuramoto Model on a Graph:}
We now consider a network of \(N\) phase oscillators (one per node of \(\mathcal{G}\)), whose phases are described by \(\theta = (\theta_1,\dots,\theta_N)^T \in \mathbb{R}^N\). 
When the graph is unweighted, the canonical Kuramoto equations are:
\begin{equation}
\label{eq:kuramoto}
\dot{\theta}_i 
\;=\;
\omega_i
\;-\;
\frac{K}{N}\,
\sum_{j:\,(i,j)\in\mathcal{E}}
\sin\bigl(\theta_i - \theta_j\bigr),
\quad
i=1,\dots,N.
\end{equation}
In matrix form, we can rewrite \eqref{eq:kuramoto} using the incidence matrix \(B\). 
Let us define the vector-valued function 
\(\sin(B^T \theta)\in\mathbb{R}^e\) for each component of \(B^T \theta\). Then \eqref{eq:kuramoto} becomes
\[
\dot{\theta}
\;=\;
\omega 
\;-\;
\frac{K}{N}\,B\,\sin\bigl(B^T \theta\bigr),
\]
which makes explicit the dependence on the graph structure via \(B\). 

We can further simplify analysis by \emph{removing} the natural frequencies \(\omega\). Specifically, let
\[
\Omega
\;=\;
\frac{1}{N}\sum_{i=1}^N \omega_i
\quad\text{and define}\quad
\widetilde{\theta}_i(t)
\;=\;
\theta_i(t)
\;-\;
\Omega\,t.
\]
Substituting into \eqref{eq:kuramoto}, we see that 
\[
\dot{\widetilde{\theta}}_i
\;=\;
\bigl(\omega_i - \Omega\bigr)
\;-\;
\frac{K}{N}\sum_{j}\sin\bigl(\widetilde{\theta}_i - \widetilde{\theta}_j\bigr).
\]
Hence, the \(\widetilde{\theta}_i\) equations have \emph{zero mean frequency}, which can simplify stability analysis. We will often assume, without loss of generality, that \(\sum_i \omega_i = 0\), or equivalently \(\Omega=0\). 

\subsection{Synchronization of Identical Coupled Oscillators}
\textbf{Result 1}: \textit{Consider the coupled oscillator model as defined earlier. If $\omega_i \neq \omega_j$ for some distinct $i, j \in\{1, \ldots, n\}$, then the oscillators cannot achieve phase synchronization.}
\vspace{8 pt}

\textbf{Proof:} We prove the lemma by contradiction. Assume that all oscillators are in phase synchrony $\theta_i(t)=\theta_j(t)$ for $t \geq 0$ and $i, j \in\{1, \ldots, n\}$. Then equating the dynamics, $\dot{\theta}_i(t)=\dot{\theta}_j(t)$, implies that $\omega_i=\omega_j$.
\vspace{8 pt}

Without loss of generality, let us now take the above dynamics equation with all angular frequencies of the system set to zero. (rotation of axes would help us achieve this)
\begin{equation*}
    \dot{\theta}=-\frac{K}{N} B \sin \left(B^T \theta\right) .
\end{equation*}

\textbf{Result 2:} \textit{If we take the unperturbed Kuramoto model defined over an arbitrarily connected graph with incidence matrix $B$ for any value of the coupling $K>0$, all trajectories will converge to the set of equilibrium solutions. In particular, the synchronized state is locally asymptotically stable. Moreover, the rate of approach to the synchronized state is no worse than $(2 K / \pi N) \lambda_2(L)$, meaning that in a local neighborhood around the equilibrium, trajectories converge at least as fast as \(\mathrm{e}^{-\alpha t}\).}
\vspace{8pt}

\textbf{Proof 1}: Consider the function $U_1(\theta)=1-r^2=\frac{4\left\|\sin \left(\frac{B^T \theta}{2}\right)\right\|^2}{N^2}$, where $r^2$ is as defined above. Taking the derivative along trajectories w.r.t time and using the fact that $\nabla_\theta U=\left(2 / N^2\right) B \sin \left(B^T \theta\right)$ leads to
$$
\dot{U}(\theta)=\nabla_\theta U \dot{\theta}=-\frac{2}{K N} \dot{\theta}^T \dot{\theta} \leq 0 .
$$
Therefore, the positive function $ U(\theta) \in [0,1]$ is a non-increasing function along the trajectories of the system. By using LaSalle's invariance principle we conclude that $U$ is a Lyapunov function for the system, and that all trajectories converge to the set where $\dot{\theta}$ is zero, i.e., the fixed point solutions.
\vspace{8pt}

\textbf{Proof 2}: We could use a different Lyapunov function similar to the approach above, infact a small angle approximation of the above Lyapunov function and consider the quadratic Lyapunov function candidate $U_2=\frac{1}{2} \sum_{i=1}^N \sum_{j=1}^N\left(\theta_i-\theta_j\right)^2=\theta^T L_c \theta$, where $L_c=N I-\mathbf{1 1}^{\mathbf{T}}$ is the Laplacian matrix of a complete graph. Note that $B^T \mathbf{1}=0$ and therefore, taking derivatives w.r.t time we get:
$$
\dot{U}_2=-\frac{K}{N} \theta^T B \sin \left(B^T \theta\right)=-\frac{K}{N} \theta^T B W(\phi) B^T \theta \leq 0
$$
We can also show that locally, the convergence is exponential with the rate determined by the smallest non-zero eigenvalue of the weighted Laplacian:
$$
\dot{U}_2 \leq-\frac{K}{N} \lambda_2\left(B W(\phi) B^T\right)\left\|\theta_{\mathbf{1}^{\perp}}\right\|^2 \leq-\frac{2 K}{\pi N} \lambda_2(L)\left\|\theta_{\mathbf{1}^{\perp}}\right\|^2
$$
using $\lambda_2\left(B W(\phi) B^T\right) \leq(2 / \pi) \lambda_2\left(B B^T\right)$ which proves the above result. 

We present two proofs because each one offers a different perspective. The trigonometric Lyapunov approach uses quantities more directly related to the classical \emph{order parameter} in Kuramoto theory, whereas the quadratic Lyapunov function is often more straightforward to manipulate for local exponential estimates.

\subsection{Existence and uniqueness of stable fixed points}
The fixed point equation can be written as:
$$
\theta=\left(B W\left(B^T \theta\right) B^T\right)^{\#} \frac{N \omega}{K}=L_W^{\#}\left(B^T \theta\right) \frac{N \omega}{K} .
$$
Using Brouwer's fixed point theorem (that states a continuous function that maps a non-empty compact, convex set $X$ into itself has at least one fixed-point), we can develop conditions which guarantee the existence (but not uniqueness) of the fixed point. If a fixed-point exists in any compact subset of $\theta \in\left(-\frac{\pi}{4}, \frac{\pi}{4}\right)$, it is stable, since this will ensure that $B^T \theta$ is between $-\frac{\pi}{2}$ and $\frac{\pi}{2}$. We therefore have to ensure that
$$
K>\frac{4}{\pi} N \max _{\left|\theta_i\right|<\frac{\pi}{4}}\left\|L_W^{\#}\left(B^T \theta\right)\right\|_{2}\|\omega\|_{2} .
$$

We are hence imposing bounds with respect to the 2-norm. So, we have:
$$
\left\|\left(B W\left(B^T \theta\right) B^T\right)^{\#}\right\|_2 \frac{N\|\omega\|_2}{K} \leq \frac{\pi}{4} \sqrt{N}
$$
Hence, a sufficient condition for synchronization of all oscillators can be determined in terms of a lower bound for $K$ :
where we used the fact that $\left\|\left(B W\left(B^T \theta\right) B\right)^{\#}\right\|_2=\frac{1}{\lambda_2\left(L_W\right)}$, and $\lambda_2$ is the algebraic connectivity of the (weighted) graph. A lower bound on the minimum value of $\lambda_2$ occurs for the minimum value of the weight which is $\frac{2}{\pi}$. As a result,
$$
K_L \geq 2 \frac{\sqrt{N}\|w\|_2}{\lambda_2(L)}
$$

\subsection{Bounds for existence of a unique fixed point:}
Consider the Kuramoto model for non-identical coupled oscillators with distinct natural frequencies \(\omega_i\). Define the threshold as shown and then:
\[
K_L
\;:=\;
2\,\frac{\sqrt{N}\,\|\omega\|_2}{\lambda_2(L)}.
\]
\begin{enumerate}
    \item For \(K \ge K_L\), there exists \emph{at least one} fixed point \(\theta^*\) with \(\bigl|\theta_i^*\bigr| < \tfrac{\pi}{4}\) or \(\bigl|\bigl(B^T\theta^*\bigr)_i\bigr| < \tfrac{\pi}{2}\).
    \item For 
    \[
    K 
    \;\ge\;
    \frac{\pi^2}{4}
    \,\frac{N\,\lambda_{\max}(L)\,\|\omega\|_2}{\lambda_2(L)^2},
    \]
    there is a \emph{unique} stable fixed point (up to adding any vector in the span of \(\mathbf{1}_N\)), and the order parameter is strictly increasing under this stronger coupling condition, implying a unique synchronized cluster in phase.
\end{enumerate}
In essence, sufficiently strong coupling ensures the existence (and, at an even higher threshold, uniqueness) of a stable phase-locked solution. 

These bounds show that if the coupling \(K\) is sufficiently large relative to the “spread” of intrinsic frequencies and the graph’s connectivity, then the system is guaranteed to have a \emph{unique}, phase-locked solution. For moderate \(K\), one can only guarantee that \emph{some} stable solution exists but not necessarily that it is unique.

\newpage
{\color{gray}\hrule}
\begin{center}
\section{Synchronization}
\textbf{A primer on Synchronization of Kuramoto Oscillatoras}
\end{center}
{\color{gray}\hrule}
\subsection{Notions of Synchronization}
\textbf{(Clarification: ${1}_n^{\perp}$ denotes perpendicular to  ${1}_n$)}
There are various notions of synchronization, described as follows:
\begin{enumerate}
    \item \textbf{Frequency synchrony:} A solution $\theta: \mathbb{R}_{\geq 0} \rightarrow \mathbb{T}^n$ is frequency synchronized if $\dot{\theta}_i(t)=\dot{\theta}_j(t)$ $\forall$ $t$ and $\forall$ $i$ and $j$.
    \item \textbf{Phase synchrony}: A solution $\theta: \mathbb{R}_{\geq 0} \rightarrow \mathbb{T}^n$ is phase synchronized if $\theta_i(t)=\theta_j(t)$ $\forall$ $t$ and $\forall$ $i$ and $j$.
Before moving to the next notion of synchronization, let us look at a few mathematical preliminaries. 
\vspace{8pt}

A torus $\mathbb{T}^n$ is the set consisting of $n$-copies of the circle. Let $\mathbb{G}$ be an undirected, weighted graph, with $\gamma$ $\in$ [0, $\pi$]. Then:
\begin{itemize}
    \item The $\operatorname{arc}$ subset $\bar{\Gamma}_{\text {arc }}(\gamma) \subset \mathbb{T}^n$ is the set of $\left(\theta_1, \ldots, \theta_n\right) \in \mathbb{T}^n$ such that there exists an arc of length $\gamma$ in $\mathbb{S}^1$ containing all angles $\theta_1, \ldots, \theta_n$. The set $\Gamma(\gamma)$ is the interior of $\bar{\Gamma}_{\text {arc }}(\gamma)$. 
    \item  The cohesive subset $\Delta^G(\gamma) \subseteq \mathbb{T}^n$ is
$$
\Delta^G(\gamma)=\left\{\theta \in \mathbb{T}^n|| \theta_i-\theta_j \mid \leq \gamma, \quad \text{for all edges }(i, j)\right\}
$$
\end{itemize}
\item \textbf{Phase cohesiveness}: A solution $\theta: \mathbb{R}_{\geq 0} \rightarrow \mathbb{T}^n$ is phase cohesive with respect to $\gamma>0$ if one of the following conditions holds $\forall$ $t$ :
\begin{itemize}
    \item $\theta(t) \in \Gamma_{\text {arc }}(\gamma)$; or
    \item $\theta(t) \in \Delta^G(\gamma)$, for the graph $G$.
\end{itemize}
\item \textbf{Asymptotic Synchronization}: This happens in cases where one of the above criterion is asymptotically achieved. For example, a solution $\theta: \mathbb{R}_{\geq 0} \rightarrow \mathbb{T}^n$ achieves phase synchronization if $\lim _{\iota \rightarrow \infty}\left|\theta_i(t)-\theta_j(t)\right|=0$.
\end{enumerate}

\subsection{Results based on above notions:}

\textbf{Res. 1} (Synchronization frequency). Consider the coupled oscillator model as defined originally with frequencies $\omega \in \mathbb{R}^n$ defined over a connected weighted undirected graph. If a solution achieves frequency synchronization, then it does so with a constant synchronization frequency equal to
$$
\omega_{\mathrm{sync}} \triangleq \frac{1}{n} \sum_{i=1}^n \omega_i=\operatorname{average}(\omega) .
$$

\textbf{Proof:} This fact is obtained by summing all equations $$
\dot{\theta}_i=\omega_i+\frac{K}{n} \sum_{j=1}^n \sin \left(\theta_j-\theta_i\right)
$$for $i \in\{1, \ldots, n\}$.

\vspace{8 pt}

\textbf{Res. 2} Consider a coupled oscillator model with frequencies $\omega \in {1}_n^{\perp}$ defined over a connected weighted undirected graph with incidence matrix $B$. The following statements hold:
\begin{enumerate}
\item (Jacobian:) the Jacobian of the coupled oscillator model at $\theta \in \mathbb{T}^n$ is
$$
J(\theta)=-B \operatorname{diag}\left(\left\{a_{i j} \cos \left(\theta_i-\theta_j\right)\right\}_{\{i, j\} \subset E}\right) B^{\top}
$$
\item (local stability:) if there exists an equilibrium $\theta^* \in \Delta^G(\gamma), \gamma<\pi / 2$, then
\begin{itemize}
\item $-J\left(\theta^*\right)$ is a Laplacian matrix; and
\item the equilibrium set $\left[\theta^*\right]$ (the rotation set for an equilibrium point obtained during rotation by a certain angle) is locally exponentially stable;
\end{itemize}
\item ((frequency synchronization:) if a solution $\theta(t)$ is phase cohesive in the sense that $\theta(t) \in \Delta^G(\gamma), \gamma<\pi / 2$, for all $t \geq 0$, then there exists a phase cohesive equilibrium $\theta^* \in \Delta^G(\gamma)$ and $\theta(t)$ achieves exponential frequency synchronization converging to $\left[\theta^*\right]$.
\end{enumerate}

\textbf{Proof:} Given $\theta \in \mathbb{T}^n$, we define the undirected graph $G_{\text {cosine }}(\theta)$ with the same nodes and edges as $G$ and with edge weights $a_{i j} \cos \left(\theta_i-\theta_j\right)$. Next, we compute
$$
\begin{aligned}
& \frac{\partial}{\partial \theta_i}\left(\omega_i-\sum_{j=1}^n a_{i j} \sin \left(\theta_i-\theta_j\right)\right)=-\sum_{j=1}^n a_{i j} \cos \left(\theta_i-\theta_j\right), \\
& \frac{\partial}{\partial \theta_j}\left(\omega_i-\sum_{k=1}^n a_{i k} \sin \left(\theta_i-\theta_k\right)\right)=a_{i j} \cos \left(\theta_i-\theta_j\right) .
\end{aligned}
$$
Therefore, the Jacobian is equal to negative of the Laplacian matrix of the graph $G_{\text {cosine}}(\theta)$. And, if $\left|\theta_i^*-\theta_j^*\right|<\pi / 2$ for all $\{i, j\} \in E$, then $\cos \left(\theta_i^*-\theta_j^*\right)>0$ for all $\{i, j\} \in E$, so that $G_{\text {cosine}}(\theta)$ has strictly non-negative weights so first part of 2. can be proved. 
\vspace{8 pt}

For the next part, we use the property that $J\left(\theta^*\right)$ is negative semidefinite with the nullspace ${1}_n$ arising from the rotational symmetry. All other eigenvectors are orthogonal to ${1}_n$ and have negative eigenvalues. Let's take a coordinate transformation matrix $Q \in \mathbb{R}^{(n-1) \times n}$ with orthonormal rows orthogonal to ${1}_n$,
$$
Q {1}_n=\mathbb{O}_{n-1} \text { and } Q Q^{\top}=I_{n-1}
$$
and we note that $Q J\left(\theta^*\right) Q^{\top}$ has negative eigenvalues. So in our original coordinate system, the zero eigenspace ${1}_n$ is exponentially stable, and hence the set $\left[\theta^*\right]$ is locally exponentially stable.
\vspace{8 pt}

For the last one, we can take a system $\dot{x}(t)=J(\theta(t)) x(t)$ and proceed. The associated undirected graph has time-varying yet strictly positive weights $a_{i j} \cos \left(\theta_i(t)-\theta_j(t)\right) \geq a_{i j} \cos (\gamma)>0$ for each $\{i, j\} \in E$. 






\subsection{Onset of Synchronization: (for non-identical oscillators)}
Here, we will first calculate the necessary critical gain $K_c$ for the onset of synchronization. As we are interested in studying the phase difference dynamics, we have 
\begin{equation*}
\begin{aligned}
\dot{\theta}_i - \dot{\theta}_j 
&= \omega_i - \omega_j 
\;+\;
\frac{K}{N} \Biggl\{ 
   -2 \sin\bigl(\theta_i - \theta_j\bigr)
   \;+\;
   \sum_{\substack{k = 1 \\ k \neq i, j}}^N 
      \Bigl[
         \sin\bigl(\theta_k - \theta_i\bigr)
         \;+\;
         \sin\bigl(\theta_j - \theta_k\bigr)
      \Bigr]
\Biggr\}.
\end{aligned}
\end{equation*}

If the oscillators are to at least asymptotically synchronize i.e. $\dot{\theta}_i-\dot{\theta}_j \rightarrow 0$ as $t \rightarrow$ $\infty \texttt{ }\forall \texttt{ }i, j=1, \ldots, N$, the R.H.S of the above equation must go to zero. 
\vspace{8 pt}

So, we derive a necessary (but not sufficient) condition for the onset of synchronization. Synchronization can only occur if our original Kuramoto equation has at least one fixed point, and hence we have:
\begin{equation*}
    \omega_j-\omega_i =\frac{K}{N}\left\{2 \sin \left(\theta_j-\theta_i\right)+\sum_{k=1}^N\left(\sin \left(\theta_k-\theta_i\right)\right.\right. \left.\left.+\sin \left(\theta_j-\theta_k\right)\right)\right\}
\end{equation*}

Without loss of generality let's assume $w_j$ $\geq$ $w_i$, then we have to maximize the following expression:
\begin{equation*}
E=2 \sin \left(\theta_j-\theta_i\right)+\sum_{k=1}^N \sin \left(\theta_k-\theta_i\right)+\sin \left(\theta_j-\theta_k\right)
\end{equation*}

Using elementary calculus, we can equate the first order derivatives wrt $i,j,k$ to zero, which gives us:
\begin{equation*}
\frac{\partial E}{\partial \theta_k} =\cos \left(\theta_k-\theta_i\right)-\cos \left(\theta_j-\theta_k\right)=0
\end{equation*}

From the last equation, we get $\theta_k=\frac{\theta_i+\theta_j}{2} \quad$ or $\quad \theta_i=\theta_j$. Since the latter gives us $E$ equals zero, we use the former condition and proceed.

\begin{equation}
\begin{aligned}
& 2 \cos \left(\theta_j-\theta_i\right)+\sum_{k=1, k \neq i, j}^N \cos \left(\frac{\theta_j-\theta_i}{2}\right)=0 \\
& \Rightarrow 2 \cos \left(\theta_j-\theta_i\right)+(N-2) \cos \left(\frac{\theta_j-\theta_i}{2}\right)=0 \\
& \Rightarrow 4 \cos ^2\left(\frac{\theta_j-\theta_i}{2}\right)-2+(N-2) \cos \left(\frac{\theta_j-\theta_i}{2}\right)=0
\end{aligned}
\end{equation}
Solving the quadratic equation gives us:
\begin{equation*}
\cos \left(\frac{\theta_j-\theta_i}{2}\right)=\frac{-(N-2)+\sqrt{(N-2)^2+32}}{8}
\end{equation*}
If this is our optimal $\theta_j - \theta_i$, say the maximum value of E is given by:
\begin{equation*}
E_{\max }=2 \sin \left(\theta_j-\theta_i\right)_{o p t}+2(N-2) \sin \left(\frac{\left(\theta_j-\theta_i\right)_{o p t}}{2}\right)
\end{equation*}
Thus, the critical coupling gain desired for onset of synchronization is:
\begin{equation*}
K_c=\frac{\left(\omega_j-\omega_i\right) N}{E_{\max }}
\end{equation*}
If the natural frequencies belong to a compact (closed, bounded) set, this becomes: 
\begin{equation*}
K_c=\frac{\left(\omega_{max}-\omega_{min}\right) N}{E_{\max }}
\end{equation*}

So, $K_c$ is simply the critical gain below which synchronization cannot occur. Now, the value for critical coupling given in [4] is:
\begin{equation*}
K_L=\frac{\left(\omega_{\max }-\omega_{\min }\right) N}{2(N-1)}
\end{equation*}

Comparing the denominators of the above bounds, we can say that $E_{max}$ equals to $2(N-1)$ is not possible in our case, which clearly is not possible as the phase differences $(\theta_m-\theta_n) \forall m, n=1, \ldots, N$ are not independent. Thus the onset of synchronization is not possible for all coupling gains $K$ satisfying $K_L \leq K<K_c$. \textbf{Hence, we have derived a stronger lower bound for the onset of synchronization}. 


\subsection{Sufficient condition for Synchronization}
Now, we will derive a sufficient condition for synchronization. The assumption in the analysis is that the initial phase of all oscillators lie within the set: 
\begin{equation*}
\mathcal{D}=\left\{\theta_i, \theta_j \in R|
| \theta_i-\theta_j \mid \leq \frac{\pi}{2}-2 \epsilon\right\}
\end{equation*}
We will find a lower bound on the coupling gain $K$ denoted by $K_{i n v}$ which makes this set positively invariant for all oscillators, i.e. $\theta_i-\theta_j \in \mathcal{D}$ at $\mathrm{t}=0 \Rightarrow \theta_i-\theta_j \in \mathcal{D} \forall t>0$. Then having phase-locked the oscillators in $\mathcal{D}$, we will show that the oscillators synchronize. We now, have:
\begin{equation}
\begin{aligned}
\dot{\theta}_i - \dot{\theta}_j 
&= 
K \Biggl\{
   \frac{\omega_i - \omega_j}{K}
   \;-\;
   \sin\bigl(\theta_i - \theta_j\bigr)
   \;+\;
   \frac{1}{N}
   \Bigl[
      \sum_{k=1}^{N} 
         \sin\bigl(\theta_i - \theta_j\bigr)
      \;+\;
      \sin\bigl(\theta_k - \theta_i\bigr)
      \;+\;
      \sin\bigl(\theta_j - \theta_k\bigr)
   \Bigr]
\Biggr\}.
\end{aligned}
\end{equation}

Rewrite the term:
\begin{equation*}
\frac{1}{N}\left(\sin \left(\theta_i-\theta_j\right)+\sin \left(\theta_k-\theta_i\right)+\sin \left(\theta_j-\theta_k\right)\right)
\end{equation*}
as
\begin{equation*}
\begin{aligned}
\frac{1}{N} \sin \left(\theta_i-\theta_j\right)\left(1-\frac{\cos \left(\theta_k-\frac{\left(\theta_i+\theta_j\right)}{2}\right)}{\cos \left(\frac{\left(\theta_i-\theta_j\right)}{2}\right)}\right) 
=\frac{1}{N} \sin \left(\theta_i-\theta_j\right) C_k
\end{aligned}
\end{equation*}
where $C_k$ $\in$ [0,1). Now, we have
\begin{equation}
    \begin{aligned}
\dot{\theta}_i-\dot{\theta}_j =K\left\{\frac{\omega_i-\omega_j}{K}-\sin \left(\theta_i-\theta_j\right)\left(1-\frac{1}{N} \sum_{k=1}^N C_k\right)\right\}
    \end{aligned}
\end{equation}

Now, let us state a result. 
\vspace{8 pt}

\textbf{Result 1} Consider the system dynamics as described by (3). Let all initial phase differences at $\mathrm{t}=0$ be contained in the compact set $\mathcal{D}=\left\{\theta_i, \theta_j|| \theta_i-\theta_j \mid \leq \frac{\pi}{2}-2 \epsilon \quad \forall i, j=\right.$ $1, \ldots, N\}$. Then there exists a coupling gain $K_{\text {inv }}>0$ such that $\left(\theta_i-\theta_j\right) \in \mathcal{D} \quad \forall t>0$.
\vspace{8 pt}

\textbf{Proof:} Consider a positive definite Lyapunov function
\begin{equation}
V=\frac{1}{2 K}\left(\theta_i-\theta_j\right)^2
\end{equation}
Taking the derivative of $V$ along the trajectories of (3) wrt time, we get:
$$
\begin{aligned}
& \dot{V}=\frac{1}{K}\left(\theta_i-\theta_j\right)\left(\dot{\theta}_i-\dot{\theta}_j\right) \\
& =\left(\theta_i-\theta_j\right)\left(\frac{\omega_i-\omega_j}{K}-\sin \left(\theta_i-\theta_j\right)\left(1-\frac{1}{N} \sum_{k=1}^N C_k\right)\right) \\
& \leq\left|\theta_i-\theta_j\right|\left|\frac{\omega_i-\omega_j}{K}\right|-\left(\theta_i-\theta_j\right) \sin \left(\theta_i-\theta_j\right)\left(1-\sum_{k=1}^N \frac{C_k}{N}\right) \\
& \leq\left|\theta_i-\theta_j\right|\left|\frac{\omega_i-\omega_j}{K}\right|-\left(\theta_i-\theta_j\right) \sin \left(\theta_i-\theta_j\right)\left(1-\frac{N-2}{N}\right)
\end{aligned}
$$
where we use $C_k<1$ and that $C_k=0$ for $k=i, j$. Thus the derivative can be written as
$$
\dot{V} \leq\left|\theta_i-\theta_j\right|\left|\frac{\omega_i-\omega_j}{K}\right|-\left(\theta_i-\theta_j\right) \sin \left(\theta_i-\theta_j\right) \frac{2}{N}
$$
Hence, if $K>\frac{N\left|\omega_i-\omega_j\right|}{2 \cos (2 \epsilon)}$ the derivative of Lyapunov function is negative at $\left|\theta_i-\theta_j\right|=\frac{\pi}{2}-2 \epsilon$ and thus the phase difference cannot leave the set $\mathcal{D}$. And to conclude, 
if $K=K_{\text {inv}}>$ $\frac{N\left|\omega_{\max }-\omega_{\min }\right|}{2 \cos (2\epsilon)}$, $\theta_i-\theta_j \quad \forall i=1,2, \ldots, N$ are positively invariant with respect to the compact set $\mathcal{D}$.
\vspace{8 pt}

Now, having trapped the phase differences within the desired compact set $\mathcal{D}$ by choosing a desired coupling gain, we will mathematically prove the synchronization. 
\vspace{8 pt}

\textbf{Result 2} Consider the system dynamics as described by (3). Let all initial phase differences at $\mathrm{t}=0$ be contained in the compact set $\mathcal{D}$. If the coupling gain $K$ is chosen such that $K=K_{i n v}$, then all the oscillators asymptotically synchronize i.e. $\dot{\theta}_i-\dot{\theta}_j \rightarrow$ 0 as $t \rightarrow \infty \quad \forall i, j=1, \ldots, N$
\vspace{8 pt}

\textbf{Proof:} Consider the positive function,
$$
S=\frac{1}{2} \dot{\theta}^T \dot{\theta}
$$
where $\dot{\theta}=\left[\dot{\theta}_1 \ldots \dot{\theta}_N\right]^T$ Taking the derivative of $V$ along the trajectories of (3) wrt time, we get:
$$
\begin{aligned}
& \dot{S}=\dot{\theta}_1 \ddot{\theta}_1+\dot{\theta}_2 \ddot{\theta}_2+\ldots+\dot{\theta}_n \ddot{\theta}_n \\
& =\frac{\dot{\theta}_1}{\beta}\left(\cos \left(\theta_1-\theta_2\right)\left(\dot{\theta}_2-\dot{\theta}_1\right)+\ldots+\cos \left(\theta_n-\theta_1\right)\left(\dot{\theta}_n-\dot{\theta}_1\right)\right) \\
& +\frac{\dot{\theta}_2}{\beta}\left(\cos \left(\theta_1-\theta_2\right)\left(\dot{\theta}_1-\dot{\theta}_2\right)+\ldots+\cos \left(\theta_n-\theta_2\right)\left(\dot{\theta}_n-\dot{\theta}_2\right)\right) \\
& \vdots \\
& +\frac{\dot{\theta}_n}{\beta}\left(\cos \left(\theta_2-\theta_n\right)\left(\dot{\theta}_2-\dot{\theta}_n\right)+\ldots+\cos \left(\theta_1-\theta_n\right)\left(\dot{\theta}_1-\dot{\theta}_n\right)\right)
\end{aligned}
$$
where $\beta=\frac{N}{K}$. On rearranging terms and simplifying we have that,
$$
\dot{S}=-\frac{K}{N} \sum_{j=1}^N \sum_{i=1}^N \cos \left(\theta_i-\theta_j\right)\left(\dot{\theta}_i-\dot{\theta}_j\right)^2
$$
Due to Result 1, we have that $\left(\theta_i-\theta_j\right) \in \mathcal{D}, \forall i, j$. This gives us that $\cos \left(\theta_i-\theta_j\right)>0 \forall i, j$ and hence $\dot{S} \leq 0$. Hence all angular frequencies are bounded. Consider the set $E=\left\{\theta_i-\theta_j, \dot{\theta}_i \in R \forall i, j \mid \dot{S}=0\right\}$. The set $E$ is characterized by all trajectories such that $\dot{\theta}_i=\dot{\theta}_j, \forall i, j$. Let $M$ be the largest invariant set contained in $E$. Using \textbf{LeSalle's Invariance} Principle, all trajectories starting in $\mathcal{D}$ converge to $M$ as $t \rightarrow \infty$. Hence, the oscillators synchronize asymptotically.

\subsection{Exponential Synchronization:}
Here, we will make use of the graph-theoretic notion of Kuramoto Oscillators from earlier. Building upon the previous theorem, 
$$
S=\frac{1}{2} \dot{\theta}^T \dot{\theta}
$$
The derivative of this function along trajectories of (1) can be written as
$$
\begin{aligned}
\dot{S} & =-\frac{K}{N} \dot{\theta}^T B \operatorname{diag}(\cos (\phi)) B^T \dot{\theta} \\
& =-\frac{K}{N} \dot{\theta}^T L_K(\mathcal{G}) \dot{\theta}
\end{aligned}
$$
The matrix $L_K(\mathcal{G})=B \operatorname{diag}(\cos (\phi)) B^T \in N \times N$ is the weighted Laplacian and is described as follows
$$
\begin{aligned}
L_W(\mathcal{G})_{i i} & =\sum_{k=1, k \neq i}^N \cos \left(\theta_k-\theta_i\right) \quad \forall i=1, \ldots, N \\
L_W(\mathcal{G})_{i j} & =-\cos \left(\theta_i-\theta_j\right) \quad \forall i, j=1, \ldots, N \quad i \neq j
\end{aligned}
$$
Clearly, if all phase differences $\phi \in \mathcal{D}$, then the weighted Laplacian matrix $L_K(\mathcal{G})$ is positive-semidefinite, and hence the previous result follows. In the next result, we extend this result by developing an exponential bound on the synchronization rate of the oscillators.
\vspace{8pt}

\textbf{Result 3}: Consider the dynamics of the system as described by (1). If the phase differences given by $\phi \in \mathcal{D}$ at $t=0$ and the coupling gain is selected such that $K=K_{i n v}$, then the oscillators synchronize exponentially at a rate no worse that $\sqrt{K \sin (2 \epsilon)}$.

\textbf{Proof}: It follows from the synchronization frequency result that:
$$
\Omega=\frac{\sum_{i=1}^N \dot{\theta}_i}{N}=\frac{\sum_{i=1}^N \omega_i}{N}
$$
We can write from the result in [2] that
\begin{equation}
\dot{\theta}=\Omega 1+\delta
\end{equation}
where `1' is the $N$-dimensional vector of ones associated with the zero eigenvalue of the weighted Laplacian $L_W(\mathcal{G}), \delta \in$ $R^n$ satisfies $\sum_{i=1}^N \delta=0\left(a s \sum_{i=1}^N \dot{\theta}_i=N \Omega\right)$.
Substituting (6) in the positive definite function $S$ as defined above, we have 
\begin{equation}
\frac{d\left(\delta^T \delta\right)}{d t}=-\frac{K}{N} \delta^T L_W(\mathcal{G}) \delta
\end{equation}
(the proof is obtained using the fact that $\Omega$ is invariant, and $$
1^T L_W(\mathcal{G})=0
$$ as 1 is an eigenvector associated with the zero eigenvalue of our weighted Laplacian matrix.)
\vspace{8pt}

We can easily see from above that $\delta$ exponentially converges to origin, now as this $\delta$ will fall to zero, we can hence say from (6) that the oscillators start moving with mean frequency of the group.

 As $\lambda_2\left(L_K(\mathcal{G})\right)$ is the smallest non-zero eigenvalue of the weighted Laplacian $\lambda_2\left(L_K(\mathcal{G})\right)$, we have from (7) that
$$
\begin{aligned}
\frac{d\left(\delta^T \delta\right)}{d t} & \leq-\frac{K}{N} \delta^T \lambda_2\left(L_W(\mathcal{G})\right) \delta \\
& \leq-\frac{K}{N} \delta^T \lambda_2\left(\operatorname{Bdiag}(\cos (\phi)) B^T\right) \delta \\
& \leq-\frac{K}{N} \delta^T \sin (2 \epsilon) \lambda_2\left(B B^T\right) \delta \\
& \leq-K \sin (2 \epsilon) \delta^T \delta
\end{aligned}
$$
as the $\min \{\cos (\phi)\}: \forall \phi \in \mathcal{D}=\cos \left(\frac{\pi}{2}-2 \epsilon\right)=\sin (2 \epsilon)$ and for an all-to-all connected topology $\lambda_2\left(B B^T\right)=N$. Thus the exponential convergence rate for synchronization is no worse that $\sqrt{K \sin (2 \epsilon)}$.

{\color{gray}\hrule}

\section{Applications and Simulations}

\subsection{Kuramoto in a Power Network}

The transmission network is described by an admittance matrix $Y \in \mathbb{C}^{n \times n}$ that is symmetric and sparse with line impedances $Z_{i j}=Z_{j i}$ for each branch $\{i, j\} \in E$. The network admittance matrix is sparse matrix with nonzero off-diagonal entries $Y_{i j}=-1 / Z_{i j}$ for each branch $\{i, j\} \in E$; the diagonal elements $Y_{i i}=-\sum_{j=1, j \neq i} Y_{i j}$ assure zero row-sums.

The static model is described by the following two concepts. Firstly, according to Kirchhoff's current law, the current injection at node $i$ is balanced by the current flows from adjacent nodes:
$$
I_i=\sum_{j=1}^n \frac{1}{Z_{i j}}\left(V_i-V_j\right)=\sum_{j=1}^n Y_{i j} V_j .
$$
Here, $I_i$ and $V_i$ are the phasor representations of the nodal current injections and nodal voltages, so that, for example, $V_i=\left|V_i\right| \mathrm{e}^{\mathrm{i} \theta_i}$ corresponds to the signal $\left|V_i\right|$ cos $\left(\omega_0 t+\theta_i\right)$. The complex power injection $S_i=V_i \cdot \bar{I}_i$ then satisfies the power balance equation
$$
S_i=V_i \cdot \sum_{j=1}^n \bar{Y}_{i j} \bar{V}_j=\sum_{j=1}^n \bar{Y}_{i j}\left|V_i\right|\left|V_j\right| e^{\mathrm{i}\left(\theta_i-\theta_j\right)}
$$
Next, for a lossless network:
$$
\underbrace{P_i}_{\text {active power injection }}=\sum_{j=1}^n \underbrace{a_{i j} \cdot \sin \left(\theta_i-\theta_j\right)}_{\text {active power flow from } i \text { to } j}, \quad i \in\{1, \ldots, n\}
$$
where $a_{i j}=\left|V_i\right|\left|V_j\right|\left|Y_{i j}\right|$ denotes the maximum power transfer over the transmission line $\{i, j\}$, and $P_i=\Re\left(S_i\right)$ is the active power injection into the network at node $i$ which is positive for generators and negative for loads. 

\vspace{8pt}
Now, let's describe a dynamical model for this network. Our assumption here is that every node is described by a first-order integrator with the following intuition: node $i$ speeds up (i.e., $\theta_i$ increases) when the power balance at node $i$ is positive, and slows down (i.e., $\theta_i$ decreases) when the power balance at node $i$ is negative. This intuition leads to a Kuramoto-like equation as follows:
$$
\boxed{\dot{\theta}_i=P_i-\sum_{j=1}^n a_{i j} \sin \left(\theta_i-\theta_j\right)}
$$

We will be simulating this for a 2-node generator-load network (all parameters used: \href{https://colab.research.google.com/drive/1_NYSeJ1Kvn24xmGAY9HeqVN-jbPhxPjc?usp=share_link}{here}), and the results are as shown:

\begin{figure}[ht!]
    \centering
    \begin{subfigure}[b]{0.45\textwidth}
        \centering
        \includegraphics[width=\textwidth]{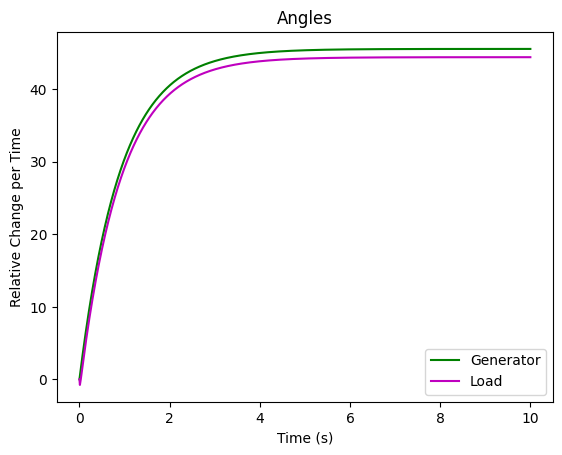}
        \label{fig:kuramoto_powernet1}
    \end{subfigure}
    \hfill
    \begin{subfigure}[b]{0.45\textwidth}
        \centering
        \includegraphics[width=\textwidth]{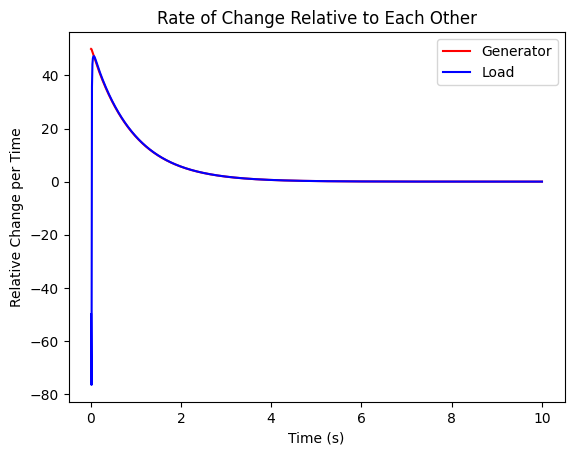}
        \label{fig:kuramoto_powernet2}
    \end{subfigure}
    \caption{ for a 2-node generator-load net}
    \label{fig:kuramoto_powernets}
\end{figure}

\subsection{Coupled Oscillator Network:}
\begin{figure}[H]
    \centering
    \includegraphics[scale = 0.7]{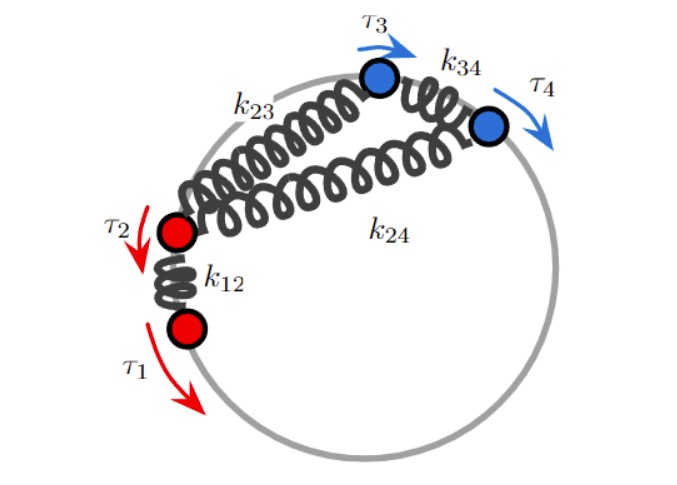}
    \caption{Springs on a ring}
    \label{fig:springs}
\end{figure}

We now analyze a system of $n$ dynamic particles constrained to rotate around a unit-radius circle and no collisions occur as shown in \ref{fig:springs}.
We assume that pairs of interacting particles $i$ and $j$ are coupled through elastic springs with stiffness $k_{i j}>0$; we set $k_{i j}=0$ if the particles are not interconnected. The elastic energy stored by the spring between particles at angles $\theta_i$ and $\theta_j$ is
$$
\begin{aligned}
\mathrm{U}_{i j}\left(\theta_i, \theta_j\right) & =\frac{k_{i j}}{2} \text {(distance)}^2=\frac{k_{i j}}{2}\left(\left(\cos \theta_i-\cos \theta_j\right)^2+\left(\sin \theta_i-\sin \theta_j\right)^2\right) \\
& =k_{i j}\left(1-\cos \left(\theta_i\right) \cos \left(\theta_j\right)-\sin \left(\theta_i\right) \sin \left(\theta_j\right)\right)=k_{i j}\left(1-\cos \left(\theta_i-\theta_j\right)\right)
\end{aligned}
$$
so that the elastic torque on particle $i$ is
$$
\mathrm{T}_i\left(\theta_i, \theta_j\right)=-\frac{\partial}{\partial \theta_i} \mathrm{U}_{i j}\left(\theta_i, \theta_j\right)=-k_{i j} \sin \left(\theta_i-\theta_j\right)
$$
From Newton's second law, we have
\begin{equation*}
    m_i \ddot{\theta}_i+d_i \dot{\theta}_i=\tau_i-\sum_{j=1}^n k_{i j} \sin \left(\theta_i-\theta_j\right),
\end{equation*}
Assuming these springs are point masses, with high damping coefficients $d$, we get
$$
\dot{\theta}_i=\omega_i-\sum_{j=1}^n a_{i j} \sin \left(\theta_i-\theta_j\right), \quad i \in\{1, \ldots, n\}
$$
with natural rotation frequencies $\omega_i=\tau_i / d$ and with coupling strengths $a_{i j}=k_{i j} / d$.
\newpage

\subsection{Vehicle Coordination: Kuramoto-Vicsek Model}
Another interesting example is the phenomenon of flocking and vehicle coordination. Let's assume that all particles have unit speed. The particle kinematics are then given by
$$
\begin{aligned}
& \dot{r}_i=e^{\mathrm{i} \theta_i}, \\
& \dot{\theta}_i=u_i(r, \theta),
\end{aligned}
$$
for $i \in\{1, \ldots, n\}$. If no control is applied, then particle $i$ travels in a straight line with orientation $\theta_i(0)$, and if $u_i=\omega_i \in \mathbb{R}$ is a nonzero constant, then particle $i$ traverses a circle with radius $1 /\left|\omega_i\right|$

The interaction among the particles is modeled by a graph $G=(\{1, \ldots, n\}, E, A)$ determined by communication and sensing patterns. Say the controllers use only relative phase information between neighboring particles (as we are mimicking biological phenomenon like the synchronization of fireflies here). Now, let's see how we can adopt potential gradient control strategies (i.e., a negative gradient flow) to coordinate the relative heading angles $\theta_i(t)-\theta_j(t)$. Let's consider a quadratic elastic spring potential to the circle $\mathrm{U}_{i j}: \mathbb{S}^1 \times \mathbb{S}^1 \rightarrow \mathbb{R}$ defined by
$$
\mathrm{U}_{i j}\left(\theta_i, \theta_j\right)=a_{i j}\left(1-\cos \left(\theta_i-\theta_j\right)\right),
$$
We can derive the affine gradient control law as follows:
$$
\dot{\theta}_i=\omega_0-K \frac{\partial}{\partial \theta_i} \sum_{\{i, j\} \subset E} \mathrm{U}_{i j}\left(\theta_i-\theta_j\right)=\omega_0-K \sum_{j=1}^n a_{i j} \sin \left(\theta_i-\theta_j\right), \quad i \in\{1, \ldots, n\} .
$$
to synchronize the heading angles of the particles for $K>0$ (gradient descent), respectively, to disperse the heading angles for $K<0$ (gradient ascent). The controlled phase dynamics above mimic animal flocking behavior. Inspired by these biological phenomena, an area of research is to study these systems in the context of tracking/flocking in swarms of autonomous vehicles. 
\vspace{8 pt}

\textbf{Simulations:} The animated result or the demonstration for flocking behaviour is \href{https://drive.google.com/file/d/1Ros44xBbO9RnbIo4VM_cFpRiPN4tbKXs/view?usp=share_link}{here}. 

\subsection{Order Parameter simulations:}
In the setup, we will consider a network of 100 oscillators with all-to-all connectivity. Shown below are the plots of order parameters for coupling constants of 0.5, 1, 2 and 3: 

\begin{figure}[H]
    \centering
    \begin{subfigure}[b]{0.45\textwidth}
        \centering
        \includegraphics[scale=0.5]{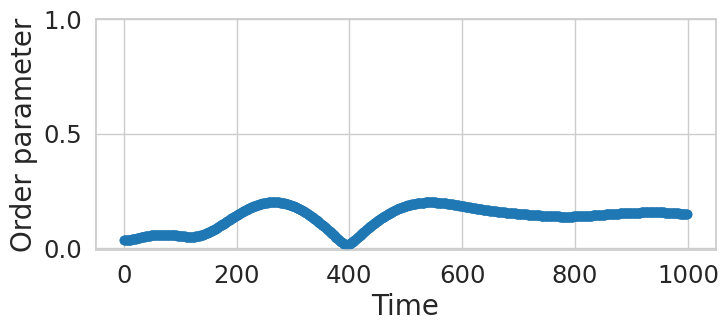}
        \caption{For K = 0.5}
        \label{fig:k1}
    \end{subfigure}
    \hfill
    \begin{subfigure}[b]{0.45\textwidth}
        \centering
        \includegraphics[scale=0.5]{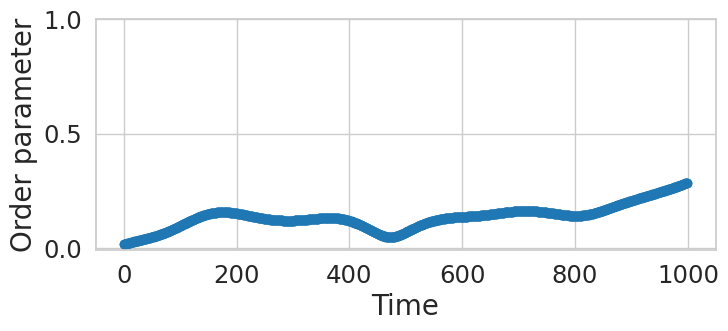}
        \caption{For K = 1}
        \label{fig:k2}
    \end{subfigure}
    \\
    \begin{subfigure}[b]{0.45\textwidth}
        \centering
        \includegraphics[scale=0.5]{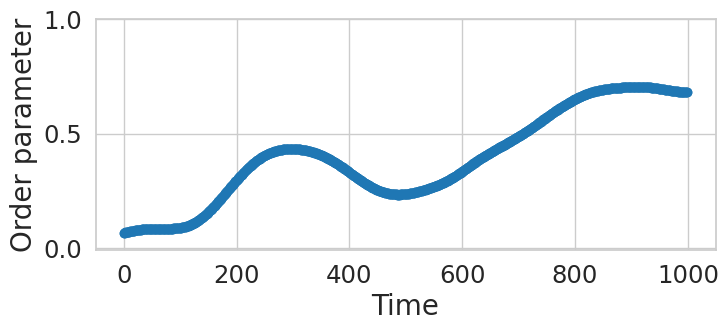}
        \caption{For K = 2}
        \label{fig:k3}
    \end{subfigure}
    \hfill
    \begin{subfigure}[b]{0.45\textwidth}
        \centering
        \includegraphics[scale=0.5]{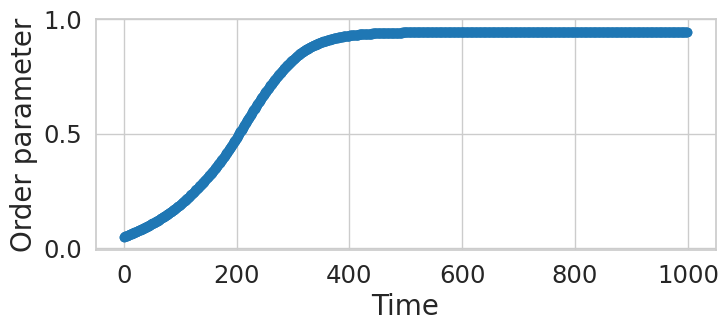}
        \caption{For K = 3}
        \label{fig:k4}
    \end{subfigure}
    \caption{Evolution of $r(t)$ with time}
    \label{fig:mainfigure}
\end{figure}








\begin{figure}[H]
    \centering
    \includegraphics[scale = 0.8]{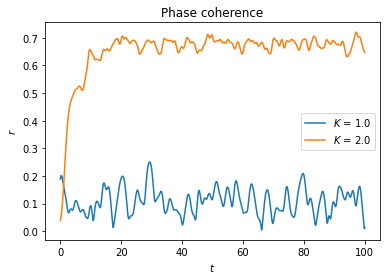}
    \caption{Showing Phase Coherence for increasing $K$}
    \label{fig:phase_coh}
\end{figure}

So, we can see that at 0.5, the order parameter is quite low and also hits zero at some point. For coupling constant of 1, the order parameter is still low but has some increasing behaviour after it hits zero. For 2, there is partial synchronization and the order parameter increases in later time intervals. For the case with $K$ as 3, we have an order parameter that gets close to 1 after some time intervals so we have full synchronization. 
\vspace{8pt}

In this case, the time series plot illustrates the same where all the trajectories eventually seem to converge. 
\begin{figure}[H]
    \includegraphics[scale = 0.6]{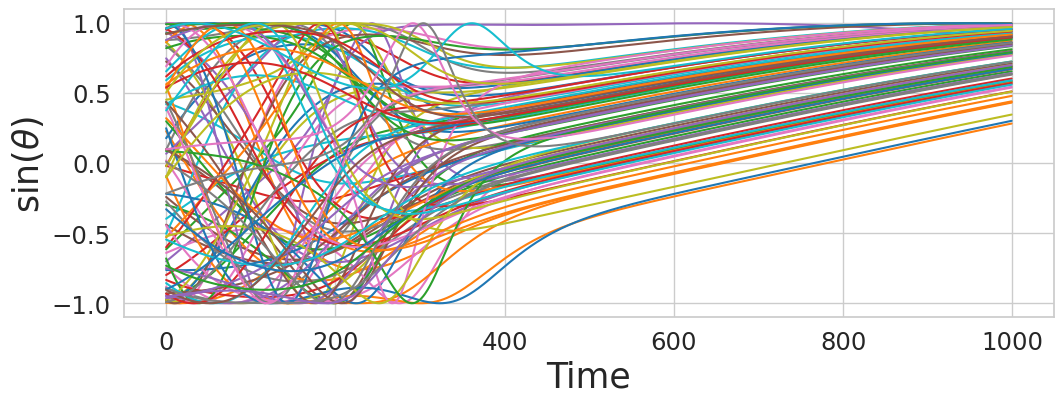}
    \caption{Evolution of trajectories}
    \label{fig:evol_traj}
\end{figure}

Now, looking plotting all the oscillators in the complex plane at different times, we get:
\begin{figure}[H]
    \includegraphics[scale = 0.5]{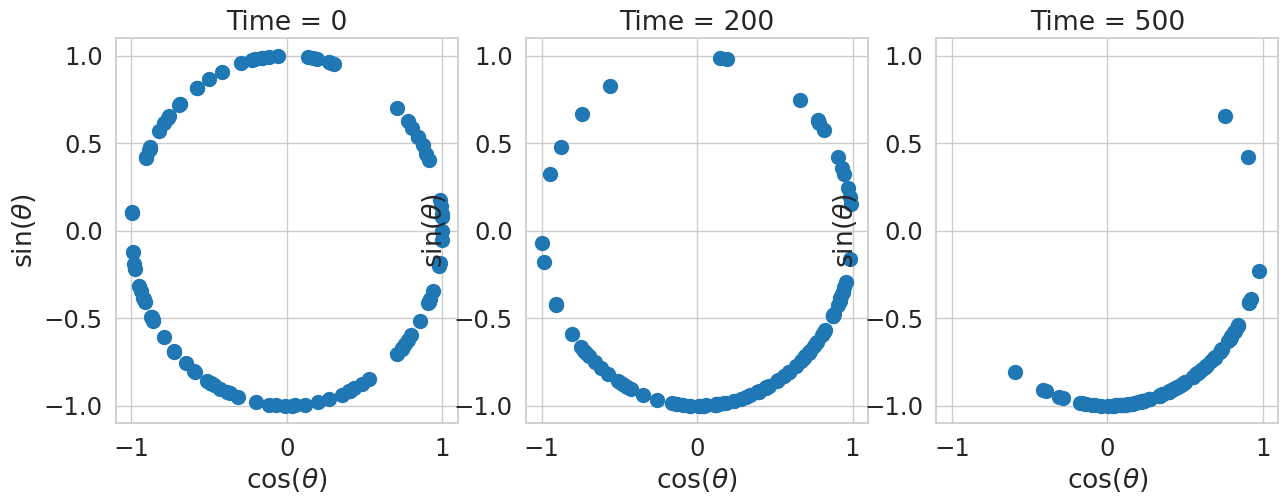}
    \caption{How the oscillators come together with time}
    \label{fig:oscil}
\end{figure}

\subsection{A Slider tool}
We also built an interactive slider tool simulated in MATLAB below: (inspiration from Cleve Moler)
\begin{figure}[H]
    \centering
    \includegraphics[scale = 0.7]{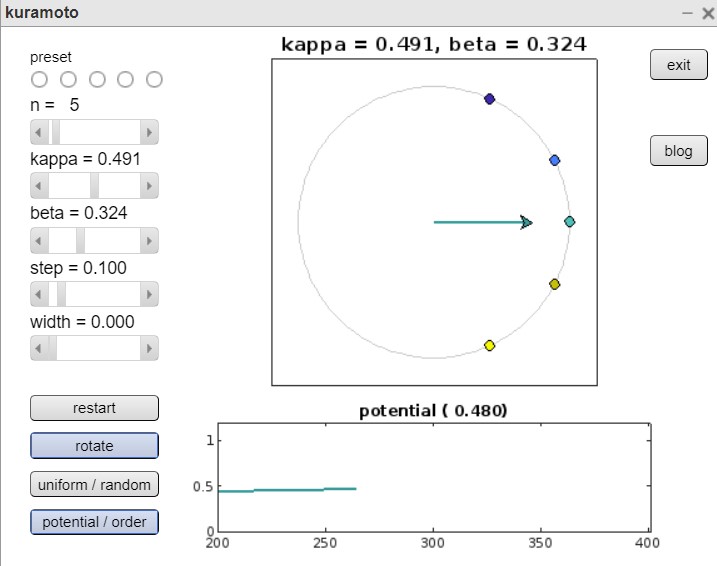}
    \caption{To plug and play}
    \label{fig:slidertool}
\end{figure}

[Code for the slider tool: \href{https://drive.google.com/file/d/1G_mHuI7GYkHlJBvvHyrN3gvxsfKoFlZb/view?usp=share_link}{here}]

\subsection{Manim Animations}
In this section, we first start out by simulating the dynamics for a given set of initial phases, and adjacency matrix $A$. Consider $n$ = 5.  

\subsubsection{Case A: When A is zero}
\begin{verbatim}
x0 = [0, 0.4, 0.8, 1.2, 1.6], w = np.array([1, 2, 3, 4, 5]), A = np.zeros((5, 5))
\end{verbatim}
We can clearly see in figure that there is no coupling. We expect the system angles to increase linearly with time at an angular velocity $w_i$.
\begin{figure}[H]
    \centering
    \includegraphics[scale = 0.7]{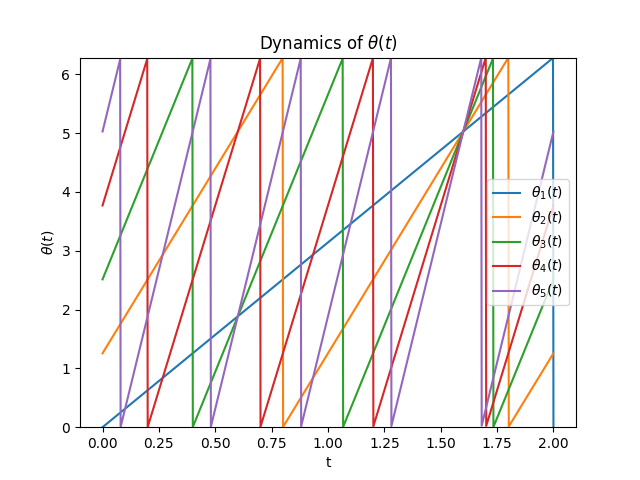}
    \caption{For $A$ = ${O_{5 X 5}}$}
    \label{fig:manim1}
\end{figure}

\subsubsection{Case B: for a Cycle Graph}
\begin{verbatim}
x0 = [0, 2 * np.pi/5 + 1/2, 4 * np.pi/5, 6 * np.pi/5 + 1/10, 8 * np.pi/5], 
w = np.array([5, 5, 5, 5, 5])
A = np.array([
    [0, 3, 0, 0, 3],
    [3, 0, 3, 0, 0],
    [0, 3, 0, 3, 0],
    [0, 0, 3, 0, 3],
    [3, 0, 0, 3, 0],])
\end{verbatim}
In this case, there is a cyclic coupling structure. Our $w$ values are constant, they are all 5. Now, let us think about the steady state behaviour of the system. 
Consider an equally spaced cyclic state. Here, $\sum_j \sin \left(\theta_i-\theta_j\right)=0$. If all the oscillators were equally space around the circle, then the corresponding all the sine terms will cancel out, and only the $w$ term will remain. So, the stable state angular velocity is just $w=5 \text{ } \forall \text{ } i$.
\vspace{8 pt}

Initially, the system is decoupled. But, in the steady state, the system becomes coupled and all the oscillators are equal spaced around the circle, and oscillate with angular velocity 5.

\begin{figure}[H]
    \centering
    \includegraphics[scale = 0.7]{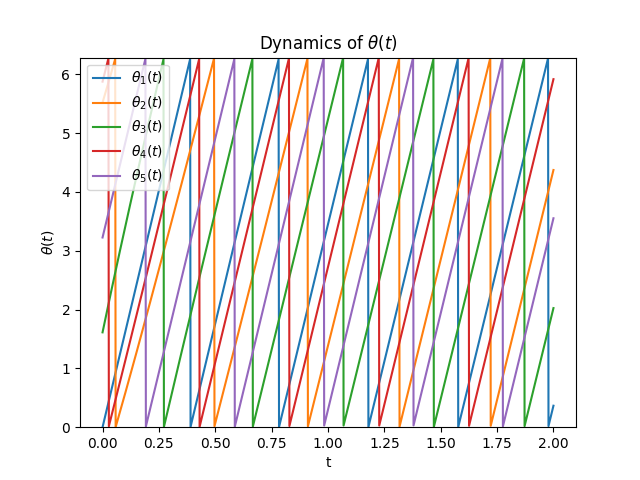}
    \caption{For a cycle graph}
    \label{fig:manim2}
\end{figure}

\subsubsection{Case C: for a Line Graph}
\begin{verbatim}
x0 = [0, 0.4, 0.8, 1.2, 1.6], w = np.array([1, 1, 1, 2.5, 2.5])
A = np.array([
    [0, 10, 0, 0, 0],
    [10, 0, 10, 0, 0],
    [0, 10, 0, 1, 0],
    [0, 0, 1, 0, 10],
    [0, 0, 0, 10, 0],])
\end{verbatim}
As the first three oscillators have the same frequency, they get coupled, while the fourth and fifth oscillators also get coupled. This phenomenon of coupling can be observed graphically also.
\begin{figure}[H]
    \centering
    \includegraphics[scale = 0.52]{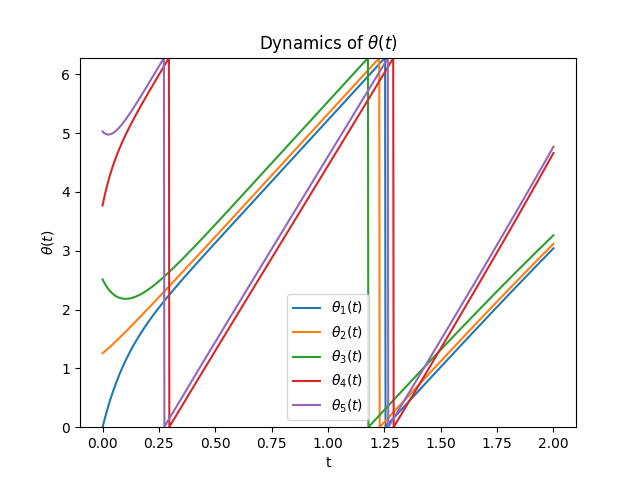}
    \caption{For a line graph}
    \label{fig:manim3}
\end{figure}
\vspace{8pt}
\newpage

[As a sidenote for cases B and C, we also tried animating the dynamics using Manim. The results can be found in \href{https://drive.google.com/drive/folders/1Qkwqj-f8JarNV3_NLYqQHc6_jg-XFZvk?usp=share_link}{here}.]




\section{Conclusions}
To conclude, we first start off with an introduction to Kuramoto Oscillators and it's broader applications. Then, we look at a graph theoretic formulation for the same, followed by a detailed discussion of various criterion for stability and synchronization of Kuramoto Oscillators. After that, we broadly analyze and experiment with various physical systems that tend to behave like Kuramoto oscillators followed by further ablation studies. \\

\footnotesize{This work was done as a part of EE6415 Non-Linear Controls project, Spring 2023.}

\section{References}
\begin{enumerate}
   \item  Steven H. Strogatz, "From Kuramoto to Crawford: exploring the onset of synchronization in populations of coupled oscillators", Physica D: Nonlinear Phenomena, Volume 143, Issues 1–4, 2000.
     \item C. Godsil and G. Royle, \emph{Algebraic Graph Theory}, Springer, 2001;  
  \item N. Biggs, \emph{Algebraic Graph Theory}, Cambridge University Press, 1993.
    \item Ali Jadbabaie, Nader Motee, Mauricio Barahona, "On the stability of the Kuramoto model of coupled nonlinear oscillators", 2005. 
    \item N. Chopra and M. W. Spong, "On Synchronization of Kuramoto Oscillators," Proceedings of the 44th IEEE Conference on Decision and Control, Seville, Spain, 2005, pp. 3916-3922, doi: 10.1109/CDC.2005.1582773.
    \item Florian Dörfler, Francesco Bullo, Synchronization in complex networks of phase oscillators: A survey, Automatica, Volume 50, Issue 6, 2014.
    \item
F. Dörfler, M. Chertkov, F. Bullo, Synchronization in complex oscillator networks and smart grids
Proc. Natl. Acad. Sci., 110 (6) (2013), pp. 2005-2010
    \item F.Bullo, Lectures on Networks and Systems, 2022, url: https://fbullo.github.io/lns.
    \item A.A. Chepizhko, V.L. Kulinskii, "On the relation between Vicsek and Kuramoto models of spontaneous synchronization", 2010. 
    \item Cooray, Gerald. The Kuramoto Model. 2008.
    \item Dirk Brockman and Steven Strogatz, "Ride my Kuramotocycle", 
    \item Cleve Moler, "Experiments With Kuramoto Oscillators", 2019.
    \item S. Strogatz. Sync: The Emerging Science of Spontaneous Order (Hyperion, New York, 2003).
    \item Balanov et al. Synchronization: From Simple to Complex. Springer-Verlag, 2009.
    \item J.A. Acebr, L.L. Bonilla, C.J.P. Vicente, F. Ritort, and R. Spigler, The Kuramoto model: A simple paradigm for synchronization phenomena, Rev. Modern Phys., 77, 137–185, 2005.
    \item Kuramoto, Yoshiki (1975). H. Araki (ed.). Lecture Notes in Physics, International Symposium on Mathematical Problems in Theoretical Physics. Vol. 39. Springer-Verlag, New York. p. 420
    \item C Huygens Oeuvres Complètes de Christiaan Huygens (Martinus Nijhoff, The Hague, The Netherlands, 1893).
optimal,” (1994)

\end{enumerate}
\end{document}